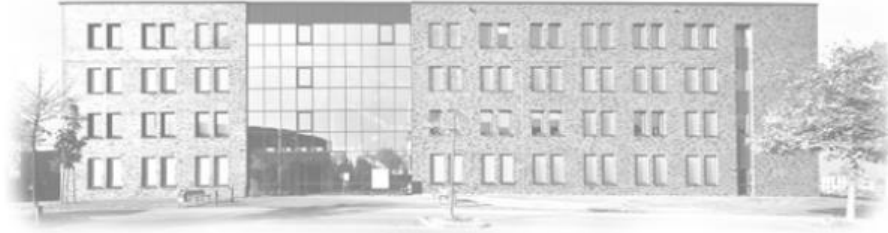

# Properties of rotation symmetric multiple-valued functions and their Reed-Muller-Fourier Spectra

**Claudio Moraga**

Lehrstuhl Informatik I

Logik in der Informatik

**Radomir S. Stanković**

Mathematical Institute

Serbian Academy of Sciences and Art

# Properties of Rotation Symmetric Multiple-Valued Functions and their Reed-Muller-Fourier Spectra

**Abstract.** The concept of rotation symmetric functions from the Boolean domain is extended to the multiple-valued (MV) domain. It is shown that symmetric functions are a subset of the rotation symmetric functions. Functions exhibiting these kinds of symmetry may be given a compact value vector representation. It is shown that the Reed-Muller-Fourier spectrum of a function preserves the kind of symmetry and therefore it may be given a compact vector representation of the same length as the compact value vector of the corresponding function. A method is presented for calculating the RMF spectrum of symmetric and rotation symmetric functions from their compact representations.

## I Introduction

As stated in [20], if a switching function has a special property, it is often realized by using few elements. At the same time, the analysis of the network realizing such a function is simplified.

Examples of desirable properties are self-duality, decomposability, linearity, maximal non-linearity, monotonicity, unateness, etc. Important classes of functions with specific properties are threshold functions, majority functions, and symmetric functions that are subject of considerations in the present Research Report. It is important to notice that many functions encountered in practice express symmetry properties. For example, functions appearing in arithmetic circuits often have symmetries [20]. Due to that, such functions can be realized with a reduced number of elements.

Symmetric functions are defined as functions that do not change their values under all possible permutations of their variables. If the permutations are restricted to certain subsets of variables, we speak about partially symmetric functions. A binary elementary symmetric function is a function $S_i^n$ such that it takes the value 1 iff exactly $i$ out of $n$ inputs are equal to 1. The set of all elementary symmetric functions $S_i^n$ for $i = 0, 1, \ldots, n$, forms a basis in terms of which an arbitrary n-variable symmetric function can be expressed as a logic or sum [20]. Elementary symmetric functions, also called simple symmetric functions, appear as voter functions in fault-tolerant computers [5].

Symmetric functions in both binary and multiple-valued cases, have been a subject of extensive study by many authors. We point out here just a grasp of references to early work in this area that we found interesting without pretending to give a deep insight into the relevant literature [1], [2], [5], [6], [11], [12], [15], [16], [17], [29].

In a study of symmetric functions, instead of restrictions to subsets of variables, we can put requirements on the allowed permutations over either subsets or all the variables. In particular, rotation symmetric functions are defined as symmetric functions under a condition that the cyclic shift of variables is a characterizing feature. They were introduced by J. Pieprzyk and C. X. Qu in the context of hashing algorithms for cryptographic applications in 1993 and 1999 [18], [19]. No applications are yet known in the multiple-valued case. These functions are however interesting on their own, as a new class of functions containing the symmetric functions. Notice that in the binary case, some rotation symmetric functions are bent functions. It has been shown by a computer search that in the binary case there are 8, 48, and 15,104, rotation symmetric bent functions in 4, 6, and 8 variables, respectively [3].

## II  BACKGROUND THEORY

The mathematical frame for this paper is the ring $(\mathbb{Z}_p, \oplus, \cdot)$. The operations are done modulo $p$, where $p$ is an integer larger than 2, not necessarily prime. The functions to be considered are mappings $f : (\mathbb{Z}_p)^n \to \mathbb{Z}_p$. For space reasons only examples for $p = 3$ and $p = 4$ will be given, but the results have general validity.

**Definition 1** *A multiple-valued function f is called symmetric, if its value does not change with any permutation of the value assignment to its arguments.*

**Definition 2** *A multiple-valued function f with more than two arguments is called **rotation symmetric** if its value does not change with any cyclic shift of patterns of value assignment to its arguments.*

**Example 1** Let $f$ be a ternary function on three arguments. Let $v_1, v_2 \in \{0, 1, 2\}$, with $v_1 \neq v_2$, be possible value assignments to the arguments.
Consider the following options:
i)   $f(v_1\, v_1\, v_1)$ has no other restriction than preserving the range,
ii)  $f(v_1\, v_1\, v_2) = f(v_1\, v_2\, v_1) = f(v_2\, v_1\, v_1)$
iii) $f(012) = f(120) = f(201)$  and
     $f(021) = f(210) = f(102)$

If $f$ satisfies the options and, moreover, $f(012) = f(021)$, then $f$ is symmetric.

If $f$ satisfies the options and, moreover, $f(012) \neq f(021)$, then $f$ is rotation symmetric: The permutations of the value assignment comprising all three values are partitioned into two disjoint cycles: (012, 120, 201) and (021, 210, 102).

**Remark 1** Symmetric functions are clearly a subset of rotation symmetric functions.

**Remark 2** When $p = 3$ and $n = 3$, except for the value assignment of type $(v_1\, v_1\, v_1)$, all permutations of the value assignment to the arguments (with two different values) are cycles of length 3 and, in the case of value assignments with three different values, they may be partitioned into two cycles of length 3, based on (012) and (021), respectively

When $p = 3$ and $n > 3$ the situation is quite different. For example if $n = 4$, let $v_1, v_2, v_3 \in \{0, 1, 2\}$, with $v_1 \neq v_2 \neq v_3$. Permutations on a value assignment $(v_1\, v_1\, v_1\, v_2)$ generate a cycle of length 4. A cyclic permutation of a value assignment $(v_1\, v_1\, v_2\, v_2)$ generates a cycle of length 4, while a cyclic permutation of a value assignment $(v_1\, v_2\, v_1\, v_2)$ generates a cycle of length 2. Value assignments comprising 1 $v_1$, 2 $v_2$'s and 1 $v_3$ generate three disjoint cycles of length 4, based on cyclic shifts of the assignments $(v_1\, v_2\, v_2\, v_3)$, $(v_2\, v_1\, v_2\, v_3)$ and $(v_2\, v_2\, v_1\, v_3)$.

**Definition 3** *For a given p and n > 2, the lexicographic first value assignment in a cycle will be its "representative".*

Since the representatives depend on the value assignments comprised in the cycles, if the permutation of a value assignment may be decomposed into disjoint cycles, the corresponding resulting representatives will automatically distinguish them.

**Definition 4** *For a given p and n > 2, if the representatives are ordered lexicographically, their position in the list – (starting with 0)– will be called their rank.*

**Example 2** *The maps of Fig. 1 show the distribution of the rank of the representatives of the vectors of value assignments when p=3 and n=3, and the values of a function which may be only symmetric or rotation symmetric depending on the ternary parameters $\alpha = \beta$ or $\alpha \neq \beta$, respectively.*

| $x_1$ | $x_2x_3$ | | | | | | | | |
|---|---|---|---|---|---|---|---|---|---|
| | 00 | 01 | 02 | 10 | 11 | 12 | 20 | 21 | 22 |
| 0 | 0 | 1 | 2 | 1 | 3 | 4 | 2 | 5 | 6 |
| 1 | 1 | 3 | 5 | 3 | 7 | 8 | 4 | 8 | 9 |
| 2 | 2 | 4 | 6 | 5 | 8 | 9 | 6 | 9 | 10 |

| $x_1$ | $x_2x_3$ | | | | | | | | |
|---|---|---|---|---|---|---|---|---|---|
| | 00 | 01 | 02 | 10 | 11 | 12 | 20 | 21 | 22 |
| 0 | 0 | 1 | 2 | 1 | 0 | $\alpha$ | 2 | $\beta$ | 1 |
| 1 | 1 | 0 | $\beta$ | 0 | 2 | 1 | $\alpha$ | 1 | 0 |
| 2 | 2 | $\alpha$ | 1 | $\beta$ | 1 | 0 | 1 | 0 | 1 |

Fig. 1: (Top) Rank of the representatives of the value assignments to the 3 ternary variables. (Bottom) A ternary function which is symmetric if $\alpha = \beta$ or rotation symmetric if $\alpha \neq \beta$.

As shown already in the early times of binary symmetric functions [22], [10], these functions may be given a compact representation as a vector of $n+1$ values ordered with increasing value of the Hamming weight of the value assignments. In the multiple-valued case a similar situation is also possible, but instead of the Hamming weight, the rank of the representatives of the value assignments in each cycle must be considered. It has been shown [16], [25], [30], that in the case of symmetric multiple-valued functions, the compact representation is given by a value vector of length κ, where

$$\kappa = \frac{(n+p-1)!}{(p-1)!n!} \quad (1)$$

Accordingly, for the ternary case, if $n = 3$, symmetric functions may be given a compact representation with a vector of length 10 and, if the function is rotation symmetric, the length of the compact vector increases to 11. There is only one permutation of value assignments that may be decomposed as a product of two disjoint cycles of length 3: (012, 120, 201) and (021, 210, 102). In summary, there are $3^{10} = 59,049$ ternary symmetric and $3^{11} = 177,147$ ternary rotation symmetric functions of three variables.. If $p = 3$ and $n = 4$ and the function is symmetric, its compact representation has a length of 15, and if the function is rotation symmetric the compact representation may reach a length of 24, if all permutations of value assignments are decomposed as products of disjoint cycles and the function takes a different value in each disjoint cycle. Clearly 24 is however still smaller than $3^4 = 81$, which would be the length of the whole list of function values. If $p = 4$ and $n = 3$ and the functions is symmetric, its compact representation has a length of 20, and if the function is rotation symmetric the compact representation may reach a length of 24. The compact representations allow to find the total number of corresponding functions. There are $3^{15} = 14,348,907$ ternary symmetric functions of 4 variables, while the number of rotation symmetric functions is $3^{24} = 282,429,536,481$.

**Example 3** *The compact representation of the ternary function of example 2 is given by*

| Representative | 000 | 001 | 002 | 011 | 012 | 021 | 022 | 111 | 112 | 122 | 222 |
|---|---|---|---|---|---|---|---|---|---|---|---|
| Rank | 0 | 1 | 2 | 3 | 4 | 5 | 6 | 7 | 8 | 9 | 10 |
| **F** | 0 | 1 | 2 | 0 | $\alpha$ | $\beta$ | 1 | 2 | 1 | 0 | 1 |

Let Σ denote the set of representatives of the vectors of value assignments to the arguments of functions for a given $p$ and $n$.

Following [25], but allowing to include rotation symmetric functions, we introduce elementary (rotation) symmetric multiple-valued functions as follows:

**Definition 5** *An elementary (rotation) symmetric function takes the value 1 at all places with value assignments belonging to the cycle of a representative from Σ.*

Table 1: Full length value vectors of elementary rotation symmetric functions for $p = 3$ and $n = 3$

| $x_1x_2x_3$ | $\phi_0$ | $\phi_1$ | $\phi_2$ | $\phi_3$ | $\phi_4$ | $\phi_5$ | $\phi_6$ | $\phi_7$ | $\phi_8$ | $\phi_9$ | $\phi_{10}$ | $x_1x_2x_3$ | $\phi_0$ | $\phi_1$ | $\phi_2$ | $\phi_3$ | $\phi_4$ | $\phi_5$ | $\phi_6$ | $\phi_7$ | $\phi_8$ | $\phi_9$ | $\phi_{10}$ |
|---|---|---|---|---|---|---|---|---|---|---|---|---|---|---|---|---|---|---|---|---|---|---|---|
| 000 | 1 | 0 | 0 | 0 | 0 | 0 | 0 | 0 | 0 | 0 | 0 | 112 | 0 | 0 | 0 | 0 | 0 | 0 | 0 | 0 | 1 | 0 | 0 |
| 001 | 0 | 1 | 0 | 0 | 0 | 0 | 0 | 0 | 0 | 0 | 0 | 120 | 0 | 0 | 0 | 0 | 1 | 0 | 0 | 0 | 0 | 0 | 0 |
| 002 | 0 | 0 | 1 | 0 | 0 | 0 | 0 | 0 | 0 | 0 | 0 | 121 | 0 | 0 | 0 | 0 | 0 | 0 | 0 | 1 | 0 | 0 | 0 |
| 010 | 0 | 1 | 0 | 0 | 0 | 0 | 0 | 0 | 0 | 0 | 0 | 122 | 0 | 0 | 0 | 0 | 0 | 0 | 0 | 0 | 0 | 1 | 0 |
| 011 | 0 | 0 | 0 | 1 | 0 | 0 | 0 | 0 | 0 | 0 | 0 | 200 | 0 | 0 | 1 | 0 | 0 | 0 | 0 | 0 | 0 | 0 | 0 |
| 012 | 0 | 0 | 0 | 0 | 1 | 0 | 0 | 0 | 0 | 0 | 0 | 201 | 0 | 0 | 0 | 0 | 1 | 0 | 0 | 0 | 0 | 0 | 0 |
| 020 | 0 | 0 | 1 | 0 | 0 | 0 | 0 | 0 | 0 | 0 | 0 | 202 | 0 | 0 | 0 | 0 | 0 | 1 | 0 | 0 | 0 | 0 | 0 |
| 021 | 0 | 0 | 0 | 0 | 0 | 1 | 0 | 0 | 0 | 0 | 0 | 210 | 0 | 0 | 0 | 0 | 0 | 1 | 0 | 0 | 0 | 0 | 0 |
| 022 | 0 | 0 | 0 | 0 | 0 | 0 | 1 | 0 | 0 | 0 | 0 | 211 | 0 | 0 | 0 | 0 | 0 | 0 | 0 | 1 | 0 | 0 | 0 |
| 100 | 0 | 1 | 0 | 0 | 0 | 0 | 0 | 0 | 0 | 0 | 0 | 212 | 0 | 0 | 0 | 0 | 0 | 0 | 0 | 0 | 0 | 1 | 0 |
| 101 | 0 | 0 | 0 | 1 | 0 | 0 | 0 | 0 | 0 | 0 | 0 | 220 | 0 | 0 | 0 | 0 | 0 | 0 | 1 | 0 | 0 | 0 | 0 |
| 102 | 0 | 0 | 0 | 0 | 0 | 1 | 0 | 0 | 0 | 0 | 0 | 221 | 0 | 0 | 0 | 0 | 0 | 0 | 0 | 0 | 1 | 0 | 0 |
| 110 | 0 | 0 | 0 | 1 | 0 | 0 | 0 | 0 | 0 | 0 | 0 | 222 | 0 | 0 | 0 | 0 | 0 | 0 | 0 | 0 | 0 | 0 | 1 |
| 111 | 0 | 0 | 0 | 0 | 0 | 0 | 0 | 1 | 0 | 0 | 0 | | | | | | | | | | | | |

Elementary rotation symmetric functions are consistent with their name: they are rotation symmetric, and therefore (for $p = 3$ and $n = 3$) they may be given a compact representation as binary valued vectors of length 11 (instead of 27). Furthermore it holds that $\phi_k(x_1,x_2,x_3) = 1$ iff the rank of the representative of the value assignment to the variables equals $k$.

## III SPECTRAL VIEW

Before moving to the spectral view of the elementary rotation symmetric functions, the Reed-Muller-Fourier transform should be recalled. This transform was introduced in 1992 [23] as a new generalization of the Reed-Muller transform to the multiple-valued domain. For a recent presentation, please refer to [27]. D.H. Green and I.S. Taylor introduced a generalization of the Reed-Muller transform to the new domain in 1976 [8]. Their generalization preserved several important properties of the binary Reed-Muller transform, however it lost the lower triangular structure and the fact of being its own inverse. The Reed-Muller-Fourier transform recovered these features besides preserving the other important properties. This was obtained by using the convolution product of Gibbs [7]. In recent comparative studies between the Reed-Muller-Transform matrix and the Pascal matrix, [26], [28], an equivalent compact definition was obtained. If for any $p$ $\mathbf{R}_1$ denotes the $p \times p$ matrix representing the basic Reed-Muller-Fourier transform, and $\mathbf{R}_1$ is expressed as $[r_{i,j}]$, $i,j \in \mathbb{Z}_p$, then

$$r_{i,j} = (-1)^j \binom{i}{j} \mod p. \tag{2}$$

Figure 2 shows the basic matrices $\mathbf{R}_1$ for $p = 3$ through 7.

$$\begin{bmatrix} 1 & 0 & 0 \\ 1 & 2 & 0 \\ 1 & 1 & 1 \end{bmatrix} ; \begin{bmatrix} 1 & 0 & 0 & 0 \\ 1 & 3 & 0 & 0 \\ 1 & 2 & 1 & 0 \\ 1 & 1 & 3 & 3 \end{bmatrix} ; \begin{bmatrix} 1 & 0 & 0 & 0 & 0 \\ 1 & 4 & 0 & 0 & 0 \\ 1 & 3 & 1 & 0 & 0 \\ 1 & 2 & 3 & 4 & 0 \\ 1 & 1 & 1 & 1 & 1 \end{bmatrix} ; \begin{bmatrix} 1 & 0 & 0 & 0 & 0 & 0 \\ 1 & 5 & 0 & 0 & 0 & 0 \\ 1 & 4 & 1 & 0 & 0 & 0 \\ 1 & 3 & 3 & 5 & 0 & 0 \\ 1 & 2 & 0 & 2 & 1 & 0 \\ 1 & 1 & 4 & 2 & 5 & 5 \end{bmatrix} ; \begin{bmatrix} 1 & 0 & 0 & 0 & 0 & 0 & 0 \\ 1 & 6 & 0 & 0 & 0 & 0 & 0 \\ 1 & 5 & 1 & 0 & 0 & 0 & 0 \\ 1 & 4 & 3 & 6 & 0 & 0 & 0 \\ 1 & 3 & 6 & 3 & 1 & 0 & 0 \\ 1 & 2 & 3 & 4 & 5 & 6 & 0 \\ 1 & 1 & 1 & 1 & 1 & 1 & 1 \end{bmatrix}$$

Fig. 2: Reed-Muller-Fourier basic matrices for $p = 3, 4, 5, 6,$ and 7.

Since the Reed-Muller Fourier transform matrix has a Kronecker product structure [Festschrift], matrices for higher dimensions are obtained as $n$-fold Kronecker product of the basic matrices with themselves.

$$\mathbf{R}_n = \mathbf{R}_1 \otimes \mathbf{R}_{n-1} = \mathbf{R}_{n-1} \otimes \mathbf{R}_1 = (\mathbf{R}_1)^{\otimes n} \tag{3}$$

**Theorem 1**.
Preliminaries:
Let $\mathbb{Z}_p$ be the domain of $p$–valued functions and let $f : (\mathbb{Z}_p)^2 \to \mathbb{Z}_p$, with value vector $\mathbf{F}$ of length $p^2$. Moreover let $g : (\mathbb{Z}_p)^2 \to \mathbb{Z}_p$, be such that $g(x_1, x_2) = f(x_2, x_1)$. Let the value vector of $g$ be $\mathbf{G}$. Furthermore, let $\mathbf{P}_2$ be a permutation matrix such that when applied upon $\mathbf{F}$ induces a permutation of its components according to the reordering of the arguments of the function. Hence $\mathbf{G} = \mathbf{P}_2 \cdot \mathbf{F}$.

Claim:
The RMF transform of a $p$-valued function on two variables preserves the order of the arguments:

$$\mathbf{R}_2 \cdot \mathbf{G} = \mathbf{R}_2 \cdot \mathbf{P}_2 \cdot \mathbf{F} = \mathbf{P}_2 \cdot \mathbf{R}_2 \cdot \mathbf{F} \mod p. \tag{4}$$

Proof:
Let $i, j \in (\mathbb{Z}_p)^2$, with $i = i_1 i_0$ and $j = j_1 j_0$.
Since $\mathbf{R}$ has a Kronecker product structure, then $\mathbf{R}_2 = \mathbf{R}_1 \otimes \mathbf{R}_1 \mod p$.
If $\mathbf{R}_2$ is expressed as $[r_{i,j}]$ then

$$r_{i,j} = \left((-1)^{j_1}\binom{i_1}{j_1}\right) \cdot \left((-1)^{j_0}\binom{i_0}{j_0}\right) = (-1)^{j_1+j_0} \frac{i_1! \cdot i_0!}{j_1!(i_1-j_1)! \, j_0!(i_0-j_0)!} \mod p \tag{5}$$

If $i_1$ and $i_0$ are exchanged, then

$$\text{modified } r_{i,j} = (-1)^{j_1+j_0} \frac{i_0! \cdot i_1!}{j_1!(i_0-j_1)! \, j_0!(i_1-j_0)!} \mod p \tag{6}$$

and if $j_1$ and $j_0$ are exchanged, then

$$\text{modified } r_{i,j} = (-1)^{j_0+j_1} \frac{i_1! \cdot i_0!}{j_0!(i_1-j_0)! \, j_1!(i_0-j_1)!} \mod p \tag{7}$$

Since scalar products and sums are commutative, equations (6) and (7) are equivalent. Therefore, exchanging $i_1$ and $i_0$ has the same effect as exchanging $j_1$ and $j_0$. Moreover exchanging $i_1$ and $i_0$ has the effect of exchanging (the corresponding) two rows of $\mathbf{R}_2$ and, similarly, exchanging $j_1$ and $j_0$ has the effect of exchanging (the corresponding) two columns of $\mathbf{R}_2$. Exchanging $i_1$ and $i_0$ corresponds to $\mathbf{P}_2 \cdot \mathbf{R}_2$, while exchanging $j_1$ and $j_0$ corresponds to $\mathbf{R}_2 \cdot \mathbf{P}_2$.

The assertion follows.

□

**Theorem 2**.
Let $n > k > 1$. Define $f$ and $g$ to be $n$-place $p$-valued functions with value vectors $\mathbf{F}$ and $\mathbf{G}$, respectively, such that for all value assignments to the arguments, $g$ equals $f$, but with transposed arguments $x_k$ and $x_{k+1}$. Let $\mathbf{P}_n$ be a permutation which when applied to $\mathbf{F}$ has the effect of transposing only the two selected arguments, i.e., $\mathbf{P}_n = (\mathbf{I}_{k-1} \otimes \mathbf{P}_2 \otimes \mathbf{I}_{n-k-1})$.
Then

$$\mathbf{R}_n \otimes \mathbf{P}_n \otimes \mathbf{F} = \mathbf{P}_n \otimes \mathbf{R}_n \otimes \mathbf{F} \mod p. \tag{8}$$

Proof:

Decompose $\mathbf{R}_n$ to match the structure of $\mathbf{P}_n$. I.e. $\mathbf{R}_n = \mathbf{R}_{k-1} \otimes \mathbf{R}_2 \otimes \mathbf{R}_{n-k-1}$, and apply it to both sides of the claim, taking advantage of the compatibility between Kronecker and matrix products [9]:

$$\mathbf{R}_n \otimes \mathbf{P}_n \otimes \mathbf{F} = (\mathbf{R}_{k-1} \otimes \mathbf{R}_2 \otimes \mathbf{R}_{n-k-1}) \cdot (\mathbf{I}_{k-1} \otimes \mathbf{P}_2 \otimes \mathbf{I}_{n-k-1}) \cdot \mathbf{F} = (\mathbf{R}_{k-1} \otimes \mathbf{R}_2 \mathbf{P}_2 \otimes \mathbf{R}_{n-k-1}) \cdot \mathbf{F} \mod p. \tag{9}$$

$$\mathbf{P}_n \otimes \mathbf{R}_n \otimes \mathbf{F} = (\mathbf{I}_{k-1} \otimes \mathbf{P}_2 \otimes \mathbf{I}_{n-k-1}) \cdot (\mathbf{R}_{k-1} \otimes \mathbf{R}_2 \otimes \mathbf{R}_{n-k-1}) \cdot \mathbf{F} = (\mathbf{R}_{k-1} \otimes \mathbf{P}_2 \mathbf{R}_2 \otimes \mathbf{R}_{n-k-1}) \cdot \mathbf{F} \mod p. \tag{10}$$

It is easy to see that the claim will be satisfied if and only if $\mathbf{P}_2\mathbf{R}_2 = \mathbf{R}_2\mathbf{P}_2$. Since this was proven in Theorem 1, the assertion follows.

□

**Theorem 3**.

Let $f$ and $g$ be $n$-place $p$-valued functions with value vectors $\mathbf{F}$ and $\mathbf{G}$, respectively, such that for all value assignments to the arguments, $g$ equals $f$, but with transposed arguments $x_k$ and $x_{k+1}$ and transposed arguments $x_h$ and $x_{h+1}$. ($n > k > h > 0$). If applied independently, let the corresponding transposition matrices be $\mathbf{P}_n^{(k)}$ and $\mathbf{P}_n^{(h)}$, respectively, leading to $\mathbf{G} = \mathbf{P}_n^{(k)} \cdot \mathbf{P}_n^{(h)} \cdot \mathbf{F}$. The following holds:

$$\mathbf{R}_n \cdot \mathbf{G} = \mathbf{P}_n^{(k)} \cdot \mathbf{P}_n^{(h)} \cdot \mathbf{R}_n \cdot \mathbf{F} \quad \mod p. \tag{11}$$

Proof:

Consider first one of the transpositions. Let $\mathbf{G}' = \mathbf{P}_n^{(h)} \cdot \mathbf{F} \mod p$.

Then from Theorem 1 follows that $\mathbf{R}_n \cdot \mathbf{G}' = \mathbf{R}_n \cdot \mathbf{P}_n^{(h)} \cdot \mathbf{F} = \mathbf{P}_n^{(h)} \cdot \mathbf{R}_n \cdot \mathbf{F} \mod p$.

Now let the second transposition be executed: $\mathbf{G} = \mathbf{P}_n^{(k)} \cdot \mathbf{G}'$.

Then from Theorem 1 follows that $\mathbf{R}_n \cdot \mathbf{G} = \mathbf{R}_n \cdot \mathbf{P}_n^{(k)} \cdot \mathbf{G}' = \mathbf{P}_n^{(k)} \cdot \mathbf{R}_n \cdot \mathbf{G}' = \mathbf{P}_n^{(k)} \cdot \mathbf{P}_n^{(h)} \cdot \mathbf{R}_n \cdot \mathbf{F} \mod p$.

□

**Theorem 4.**

Let $f$ and $g$ be $n$-place $p$-valued functions with value vectors $\mathbf{F}$ and $\mathbf{G}$, respectively, such that for all value assignments to the arguments, $g$ equals $f$, but with permuted arguments. Let $\mathbf{P}_n$ be a permutation matrix, which when applied to $\mathbf{F}$ has the same effect as permuting the corresponding arguments.
Then

$$\mathbf{R}_n \cdot \mathbf{G} = \mathbf{R}_n \cdot \mathbf{P}_n \cdot \mathbf{F} = \mathbf{P}_n \cdot \mathbf{R}_n \cdot \mathbf{F} \quad \mod p. \tag{12}$$

Proof:

Recall that any permutation of an ordered set of arguments may be decomposed as a product of disjoint cycles, and any (non-singleton) cycle may be obtained with a cascade of transpositions. (Furthermore, any transposition may be realized as a cascade of neighbor transpositions; e.g. (1 4) = (4 3)(3 2)(1 2)(2 3)(3 4).) Apply accordingly Theorems 2 and 3 as many times as needed.

□

**Theorem 5.**

Let $n > 2$. The RMF spectrum of an $n$-place $p$-valued rotation symmetric function is rotation symmetric.

Proof:

Recall that a $p$-valued function is rotation symmetric iff it is invariant with respect to any cyclic permutation of the value assignment of its arguments.

Let $\mathbf{F}$ be the value vector of a rotation symmetric function and let $\mathbf{P}_n$ be equivalent to a random cyclic permutation of the value assignment of its arguments. Then, since $\mathbf{F}$ is the value vector of a rotation symmetric function,

$$\mathbf{F} = \mathbf{P}_n \cdot \mathbf{F}. \tag{13}$$

From Theorem 4,

$$\mathbf{R}_n \cdot \mathbf{F} = \mathbf{R}_n \cdot \mathbf{P}_n \cdot \mathbf{F} = \mathbf{P}_n \cdot \mathbf{R}_n \cdot \mathbf{F} \quad \mod p.$$

Therefore $\mathbf{R}_n \cdot \mathbf{F} \mod p$ is rotation symmetric.

□

## III  COMPUTATIONAL ASPECTS

The Reed-Muller-Fourier spectrum of the elementary rotation symmetric functions on $n$ variables will be (first) calculated by multiplying the $p^n \times p^n$ transform matrix $\mathbf{R}_n$ with the value vector $\mathbf{\Phi}$ (of length $p^n$) of the corresponding function. In what follows, angular parentheses will be used to denote the compact representation of a vector. Taking in account Theorem 1, a rotation symmetric function will be shown in its compact representation ⟨**F**⟩. Similarly for Reed-Muller-Fourier spectra. For the full vectors, the notation $SB(n,k)$ introduced in [21] and used in [25] will be adopted, where $n$ denotes the number of variables of the function, and $k$ indicates the rank of the representative of the characterizing value assignment to the variables.

For $p = 3$ and $n = 3$, the compact ⟨$SB(3,k)$⟩ functions specifications are shown in Table 2.

Table 2: Compact RMF-spectrum ⟨$SB(3,k)$⟩ of elementary rotation symmetric functions for $p = 3$ and $n = 3$.

| $x_1x_2x_3$ (repr.) | ⟨$SB(3,0)$⟩ | ⟨$SB(3,1)$⟩ | ⟨$SB(3,2)$⟩ | ⟨$SB(3,3)$⟩ | ⟨$SB(3,4)$⟩ | ⟨$SB(3,5)$⟩ | ⟨$SB(3,6)$⟩ | ⟨$SB(3,7)$⟩ | ⟨$SB(3,8)$⟩ | ⟨$SB(3,9)$⟩ | ⟨$SB(3,10)$⟩ |
|---|---|---|---|---|---|---|---|---|---|---|---|
| 000 | 1 | 0 | 0 | 0 | 0 | 0 | 0 | 0 | 0 | 0 | 0 |
| 001 | 1 | 2 | 0 | 0 | 0 | 0 | 0 | 0 | 0 | 0 | 0 |
| 002 | 1 | 1 | 1 | 0 | 0 | 0 | 0 | 0 | 0 | 0 | 0 |
| 011 | 1 | 1 | 0 | 1 | 0 | 0 | 0 | 0 | 0 | 0 | 0 |
| 012 | 1 | 0 | 1 | 2 | 2 | 0 | 0 | 0 | 0 | 0 | 0 |
| 021 | 1 | 0 | 1 | 2 | 0 | 2 | 0 | 0 | 0 | 0 | 0 |
| 022 | 1 | 2 | 2 | 1 | 1 | 1 | 1 | 0 | 0 | 0 | 0 |
| 111 | 1 | 0 | 0 | 0 | 0 | 0 | 0 | 2 | 0 | 0 | 0 |
| 112 | 1 | 2 | 1 | 2 | 2 | 2 | 0 | 1 | 1 | 0 | 0 |
| 122 | 1 | 1 | 2 | 2 | 0 | 0 | 1 | 2 | 1 | 2 | 0 |
| 222 | 1 | 0 | 0 | 0 | 0 | 0 | 0 | 1 | 0 | 0 | 1 |

**Theorem 6**

Let $f$ be a rotation symmetric $p$-valued function specified by a value vector **F** and a compact value vector ⟨**F**⟩. The Reed-Muller-Fourier compact spectrum ⟨**S**$_f$⟩, of $f$, may be obtained directly from ⟨**F**⟩ based on the ⟨$SB(n,k)$⟩ vectors as follows:

$$⟨\mathbf{S}_f⟩ = ⟨\mathbf{F}⟩_j \cdot ⟨SB(n,j)⟩ \bmod 3, \tag{14}$$

where ⟨**F**⟩$_j$ denotes the $j$-th component of ⟨**F**⟩.

Proof:

Let $\phi_j$ denote the $j$-th elementary rotation symmetric function and **Φ**$_j$ its value vector.

The Reed-Muller-Fourier spectrum of $f$ with value vector **F** is given by

$$\begin{aligned}
\mathbf{S}_f &= \mathbf{R}_n \cdot \mathbf{F} \bmod 3 \\
&= \mathbf{R}_n \sum_j ⟨\mathbf{F}⟩_j \cdot \mathbf{\Phi}_j \bmod 3 \\
&= \sum_j \mathbf{R}_n \cdot ⟨\mathbf{F}⟩_j \cdot \mathbf{\Phi}_j \bmod 3 \\
&= \sum_j ⟨\mathbf{F}⟩_j \cdot \mathbf{R}_n \cdot \mathbf{\Phi}_j \bmod 3 \\
&= \sum_j ⟨\mathbf{F}⟩_j \cdot \mathbf{S}_{\phi_j} = \sum_j ⟨\mathbf{F}⟩_j \cdot SB(n,j) \bmod 3. \tag{15}
\end{aligned}$$

Then $\forall\ w$ belonging to a value assignment cycle with rank $j$ holds:

$$\mathbf{S}_f(w) = \sum_j ⟨\mathbf{F}⟩_j \cdot \mathbf{S}_{\phi_j}(w) \bmod 3, \tag{16}$$

from where

$$⟨\mathbf{S}_f⟩ = \sum_j ⟨\mathbf{F}⟩_j \cdot ⟨SB(n,j)⟩ \bmod 3. \tag{17}$$

□

Since the Reed-Muller-Fourier transform is its own inverse, the following holds:

**Corollary 6.1**

$$⟨\mathbf{F}⟩ = \sum_j ⟨\mathbf{S}_f⟩_j \cdot ⟨SB(n,j)⟩ \bmod 3 \tag{18}$$

**Example 4** *Consider the function of Example 2, with α = 1 and β = 0. The function is clearly rotation symmetric.*

| $x_1$ | \multicolumn{9}{c}{$x_2x_3$} |
|---|---|---|---|---|---|---|---|---|---|
|  | 00 | 01 | 02 | 10 | 11 | 12 | 20 | 21 | 22 |
| 0 | 0 | 1 | 2 | 1 | 0 | 1 | 2 | 0 | 1 |
| 1 | 1 | 0 | 0 | 0 | 2 | 1 | 1 | 1 | 0 |
| 2 | 2 | 1 | 1 | 0 | 1 | 0 | 1 | 0 | 1 |

Its RMF spectrum is:

| $x_1$ | \multicolumn{9}{c}{$x_2x_3$} |
|---|---|---|---|---|---|---|---|---|---|
|  | 00 | 01 | 02 | 10 | 11 | 12 | 20 | 21 | 22 |
| 0 | 0 | 2 | 0 | 2 | 1 | 1 | 0 | 2 | 2 |
| 1 | 2 | 1 | 2 | 1 | 1 | 0 | 1 | 0 | 2 |
| 2 | 0 | 1 | 2 | 2 | 0 | 2 | 2 | 2 | 0 |

The compact representation of $\langle \mathbf{F} \rangle$ and of the compact RMF-spectra of elementary rotation symmetric functions scaled by the non-zero coefficients of $\langle \mathbf{F} \rangle$ according to Theorem 6, is the following:

| $x_1x_2x_3$ (repr.) | $\langle \mathbf{F} \rangle$ | $\langle SB(3,1) \rangle$ | $2\langle SB(3,2) \rangle$ | $\langle SB(3,4) \rangle$ | $\langle SB(3,6) \rangle$ | $2\langle SB(3,7) \rangle$ | $\langle SB(3,8) \rangle$ | $\langle SB(3,10) \rangle$ | $\oplus$ |
|---|---|---|---|---|---|---|---|---|---|
| 000 | 0 | 0 | 0 | 0 | 0 | 0 | 0 | 0 | 0 |
| 001 | 1 | 2 | 0 | 0 | 0 | 0 | 0 | 0 | 2 |
| 002 | 2 | 1 | 2 | 0 | 0 | 0 | 0 | 0 | 0 |
| 011 | 0 | 1 | 0 | 0 | 0 | 0 | 0 | 0 | 1 |
| 012 | 1 | 0 | 2 | 2 | 0 | 0 | 0 | 0 | 1 |
| 021 | 0 | 0 | 2 | 0 | 0 | 0 | 0 | 0 | 2 |
| 022 | 1 | 2 | 1 | 1 | 1 | 0 | 0 | 0 | 2 |
| 111 | 2 | 0 | 0 | 0 | 0 | 1 | 0 | 0 | 1 |
| 112 | 1 | 2 | 2 | 2 | 0 | 2 | 1 | 0 | 0 |
| 122 | 0 | 1 | 1 | 0 | 1 | 1 | 1 | 0 | 2 |
| 222 | 1 | 0 | 0 | 0 | 0 | 2 | 0 | 1 | 0 |

The columns are scaled according to the elements of $\langle \mathbf{F} \rangle$. It may be seen that the column labelled "$\oplus$" represents the results of Eq. (17), returning $\langle \mathbf{S}_f \rangle$.

Similarly, Corollary 6.1 may be illustrated as follows, where the columns are now scaled by the non-zero coefficients of $\langle \mathbf{S}_f \rangle$ and the column labelled "$\oplus$" represents the results of Eq. (18), returning $\langle \mathbf{F} \rangle$.

| $x_1x_2x_3$ (repr.) | $\langle \mathbf{S}_f \rangle$ | $2\langle SB(3,1) \rangle$ | $\langle SB(3,3) \rangle$ | $\langle SB(3,4) \rangle$ | $2\langle SB(3,5) \rangle$ | $2\langle SB(3,6) \rangle$ | $\langle SB(3,7) \rangle$ | $2\langle SB(3,9) \rangle$ | $\oplus$ |
|---|---|---|---|---|---|---|---|---|---|
| 000 | 0 | 0 | 0 | 0 | 0 | 0 | 0 | 0 | 0 |
| 001 | 2 | 1 | 0 | 0 | 0 | 0 | 0 | 0 | 1 |
| 002 | 0 | 2 | 0 | 0 | 0 | 0 | 0 | 0 | 2 |
| 011 | 1 | 2 | 1 | 0 | 0 | 0 | 0 | 0 | 0 |
| 012 | 1 | 0 | 2 | 2 | 0 | 0 | 0 | 0 | 1 |
| 021 | 2 | 0 | 2 | 0 | 1 | 0 | 0 | 0 | 0 |
| 022 | 2 | 1 | 1 | 1 | 2 | 2 | 0 | 0 | 1 |
| 111 | 1 | 0 | 0 | 0 | 0 | 0 | 2 | 0 | 2 |
| 112 | 0 | 1 | 2 | 2 | 1 | 0 | 1 | 0 | 1 |
| 122 | 2 | 2 | 2 | 0 | 0 | 2 | 2 | 1 | 0 |
| 222 | 0 | 0 | 0 | 0 | 0 | 0 | 1 | 0 | 1 |

**Remark 3** Notice that if $p = 3$ and $n = 3$, there are 59,049 symmetric functions [25] and $3^{11} = 177,147$ rotation symmetric functions. However the 11 elementary rotation symmetric functions must be calculated *only once*, to move between the functions and their Reed-Muller-Fourier spectra with Theorem 6 and its Corollary. The computational complexity of Theorem 6 is in $O(\kappa^2)$, but all weighted sums may be run in parallel.

**Lemma 1**: For a given $p$ and a given $n > 2$ holds the following:
  i) The sum mod $p$ of two symmetric functions is symmetric
  ii) The sum mod $p$ of a symmetric and a rotation symmetric functions is rotation symmetric
  iii) The sum mod $p$ of two rotation symmetric functions may be
  - rotation symmetric with the same number of cycles of common permutations with different values
  - rotation symmetric with a smaller number of cycles with different values if the entries of all cycles of some permutations add up (mod $p$) to a common value
  - just symmetric if the entries of all cycles of each permutation add up (mod p) to a common value

An example instead of a formal proof:

Let $p = 3$ and $n = 4$. It may be shown that all value assignments to the arguments may be decomposed into 24 cycles. Let the assignment $A_1$ comprise the values (0, 0, 1, 2), whose permutations may be partitioned into three disjoint cycles $c_{1,1}$, $c_{1,2}$ and $c_{1,3}$ as follows:

$c_{1,1} = (0012–0120–1200–2001)$, $c_{1,2} = (0102–1020–0201–2010)$ and $c_{1,3} = (1002–0021–0210–2100)$

Similarly, let $A_2$ comprise the values (0, 0, 2, 2), whose permutations may be partitioned into two disjoint cycles as follows:

$c_{2,1} = (0022–0220–2200–2002)$ and $c_{2,2} = (0202–2020)$

Let $F_1$, $F_2$, $F_3$, $F_4$, $F_5$ and $F_6$ be the value vectors of ternary 4-place (rotation) symmetric functions, which are equal respectively at all places, except at the places in the above cycles. Let $F_1$ and $F_2$ be straight symmetric and the other functions, rotation symmetric. Accordingly, the relevant parts of these functions may be represented as follows:

|   |   | $c_{1,1}$ | $c_{1,2}$ | $c_{1,3}$ |   | $c_{2,1}$ | $c_{2,2}$ |   |
|---|---|---|---|---|---|---|---|---|
| $\langle F_1 \rangle$ | ………. | 1 | 1 | 1 | ………. | 2 | 2 | ………. |
| $\langle F_2 \rangle$ | ………. | 0 | 0 | 0 | ………. | 1 | 1 | ………. |
| $\langle F_3 \rangle$ | ………. | 2 | 1 | 0 | ………. | 2 | 1 | ………. |
| $\langle F_4 \rangle$ | ………. | 2 | 0 | 1 | ………. | 2 | 1 | ………. |
| $\langle F_5 \rangle$ | ………. | 1 | 2 | 0 | ………. | 2 | 1 | ………. |
| $\langle F_6 \rangle$ | ………. | 0 | 2 | 1 | ………. | 2 | 0 | ………. |
| $\langle F_1 \rangle \oplus \langle F_2 \rangle$ | ………. | 1 | 1 | 1 | ………. | 0 | 0 | ………. |
| $\langle F_1 \rangle \oplus \langle F_3 \rangle$ | ………. | 0 | 2 | 1 | ………. | 1 | 0 | ………. |
| $\langle F_4 \rangle \oplus \langle F_5 \rangle$ | ………. | 0 | 2 | 1 | ………. | 1 | 2 | ………. |
| $\langle F_3 \rangle \oplus \langle F_5 \rangle$ | ………. | 0 | 0 | 0 | ………. | 1 | 2 | ………. |
| $\langle F_5 \rangle \oplus \langle F_6 \rangle$ | ………. | 1 | 1 | 1 | ………. | 1 | 1 | ………. |

It is easy to see that $\langle F_1 \rangle \oplus \langle F_2 \rangle$ produces another symmetric function; $\langle F_1 \rangle \oplus \langle F_3 \rangle$ as well as $\langle F_4 \rangle \oplus \langle F_5 \rangle$ produce new rotation symmetric functions, while $\langle F_3 \rangle \oplus \langle F_5 \rangle$ produces a "weaker" rotation symmetric function since only for the cycles of the assignments of $A_2$ takes different values. Finally, $\langle F_5 \rangle \oplus \langle F_6 \rangle$, the sum mod 3 of two rotation symmetric functions, produces a straight symmetric function.

**Remark 4** In the case of $p = 4$ and $n = 3$, (minimum size of $n$ required to have one permutation of a value assignment to be decomposable into two disjoint cycles of length 3, and therefore to have the possibility of building rotation symmetric functions), $\kappa$ equals 20. There are $4^{20} \approx 10^{12}$ symmetric functions. There are 24 elementary rotation symmetric functions, each compact value vector is of length 24. It follows that there are $4^{24} \approx 256 \cdot 10^{12}$ quaternary rotation symmetric functions of three variables.

Representatives, ranks and cycles are summarized in Table 3.

Table 3: Disjoint cycles of value assignments for $p = 4$ and $n = 3$

| repr. | rank | cycle | repr. | rank | cycle | repr. | rank | cycle |
|---|---|---|---|---|---|---|---|---|
| 000 | 0 | (000) | 022 | 8 | (022-220-202) | 122 | 16 | (122-221-212) |
| 001 | 1 | (001-010-100) | 023 | 9 | (023-230-302) | 123 | 17 | (123-231-312) |
| 002 | 2 | (002-020-200) | 031 | 10 | (031-310-103) | 132 | 18 | (132-321-213) |
| 003 | 3 | (003-030-300) | 032 | 11 | (032-320-203) | 133 | 19 | (133-331-313) |
| 011 | 4 | (011-110-101) | 033 | 12 | (033-330-303) | 222 | 20 | (222) |
| 012 | 5 | (012-120-201) | 111 | 13 | (111) | 223 | 21 | (223-232-322) |
| 013 | 6 | (013-130-301) | 112 | 14 | (112-121-211) | 233 | 22 | (233-332-323) |
| 021 | 7 | (021-210-102) | 113 | 15 | (113-131-311) | 333 | 23 | (333) |

Notice that the cycles with ranks 5 and 7; 6 and 10; 9 and 11; as well as 17 and 18, are disjoint cycles obtained from permutations of length 6 where the value assignments comprise three different values. If a function has to be rotation symmetric, then it should have different values at least in one pair of these disjoint cycles (but a constant value in each cycle).

The distribution of ranks for $p = 4$ and $n = 3$ is shown in the following map:

| | | | | | | | | | $x_2 x_3$ | | | | | | | |
|---|---|---|---|---|---|---|---|---|---|---|---|---|---|---|---|---|
| $x_1$ | 00 | 01 | 02 | 03 | 10 | 11 | 12 | 13 | 20 | 21 | 22 | 23 | 30 | 31 | 32 | 33 |
| 0 | 0 | 1 | 2 | 3 | 1 | 4 | 5 | 6 | 2 | 7 | 8 | 9 | 3 | 10 | 11 | 12 |
| 1 | 1 | 4 | 7 | 10 | 4 | 13 | 14 | 15 | 5 | 14 | 16 | 17 | 6 | 15 | 18 | 19 |
| 2 | 2 | 5 | 8 | 11 | 7 | 14 | 16 | 18 | 8 | 16 | 20 | 21 | 9 | 17 | 21 | 22 |
| 3 | 3 | 6 | 9 | 12 | 10 | 15 | 17 | 19 | 11 | 18 | 21 | 22 | 12 | 19 | 22 | 23 |

Before analyzing examples with $p = 4$ and $n = 3$, the elementary rotation symmetric functions and their RMF spectra are needed. Notice that in the compact representation, $\phi_k(x_1 x_2 x_3) = 1$ iff $\text{rank}(x_1 x_2 x_3) = k$. An explicit table is not needed.

The RMF spectra of the elementary rotation symmetric functions are needed before working the following examples. For space reasons in the next tables let $i_k$ denote $\langle SB(4,k) \rangle$.

| rep. | k | $i_1$ | $i_2$ | $i_3$ | $i_4$ | $i_5$ | $i_6$ | rep. | k | $i_1$ | $i_2$ | $i_3$ | $i_4$ | $i_5$ | $i_6$ | rep. | k | $i_1$ | $i_2$ | $i_3$ | $i_4$ | $i_5$ | $i_6$ |
|---|---|---|---|---|---|---|---|---|---|---|---|---|---|---|---|---|---|---|---|---|---|---|---|
| 000 | 0 | 0 | 0 | 0 | 0 | 0 | 0 | 022 | 8 | 0 | 2 | 0 | 0 | 2 | 0 | 122 | 16 | 3 | 2 | 0 | 0 | 1 | 0 |
| 001 | 1 | 3 | 0 | 0 | 0 | 0 | 0 | 023 | 9 | 3 | 0 | 3 | 2 | 2 | 2 | 123 | 17 | 2 | 0 | 3 | 3 | 1 | 2 |
| 002 | 2 | 2 | 1 | 0 | 0 | 0 | 0 | 031 | 10 | 0 | 3 | 3 | 3 | 0 | 0 | 132 | 18 | 2 | 0 | 3 | 3 | 2 | 1 |
| 003 | 3 | 1 | 3 | 3 | 0 | 0 | 0 | 032 | 11 | 3 | 0 | 3 | 2 | 1 | 0 | 133 | 19 | 1 | 2 | 2 | 3 | 0 | 0 |
| 011 | 4 | 2 | 0 | 0 | 1 | 0 | 0 | 033 | 12 | 2 | 2 | 2 | 1 | 3 | 3 | 222 | 20 | 2 | 3 | 0 | 0 | 2 | 0 |
| 012 | 5 | 1 | 1 | 0 | 2 | 3 | 0 | 111 | 13 | 1 | 0 | 0 | 3 | 0 | 0 | 223 | 21 | 1 | 1 | 3 | 0 | 1 | 2 |
| 013 | 6 | 0 | 3 | 3 | 3 | 1 | 1 | 112 | 14 | 0 | 1 | 0 | 1 | 3 | 0 | 233 | 22 | 0 | 3 | 2 | 1 | 2 | 1 |

| rep. | k | $i_1$ | $i_2$ | $i_3$ | $i_4$ | $i_5$ | $i_6$ | rep. | k | $i_1$ | $i_2$ | $i_3$ | $i_4$ | $i_5$ | $i_6$ | rep. | k | $i_1$ | $i_2$ | $i_3$ | $i_4$ | $i_5$ | $i_6$ |
|---|---|---|---|---|---|---|---|---|---|---|---|---|---|---|---|---|---|---|---|---|---|---|---|
| 021 | 7 | 1 | 1 | 0 | 2 | 0 | 0 | 113 | 15 | 3 | 3 | 3 | 3 | 1 | 1 | 333 | 23 | 3 | 1 | 1 | 3 | 1 | 1 |

| rep. | k | $i_7$ | $i_8$ | $i_9$ | $i_{10}$ | $i_{11}$ | $i_{12}$ | rep. | k | $i_7$ | $i_8$ | $i_9$ | $i_{10}$ | $i_{11}$ | $i_{12}$ | rep. | k | $i_7$ | $i_8$ | $i_9$ | $i_{10}$ | $i_{11}$ | $i_{12}$ |
|---|---|---|---|---|---|---|---|---|---|---|---|---|---|---|---|---|---|---|---|---|---|---|---|
| 000 | 0 | 0 | 0 | 0 | 0 | 0 | 0 | 022 | 8 | 2 | 1 | 0 | 0 | 0 | 0 | 122 | 16 | 1 | 1 | 0 | 0 | 0 | 0 |
| 001 | 1 | 0 | 0 | 0 | 0 | 0 | 0 | 023 | 9 | 1 | 3 | 3 | 0 | 0 | 0 | 123 | 17 | 2 | 3 | 3 | 1 | 0 | 0 |
| 002 | 2 | 0 | 0 | 0 | 0 | 0 | 0 | 031 | 10 | 1 | 0 | 0 | 1 | 0 | 0 | 132 | 18 | 1 | 3 | 0 | 2 | 3 | 0 |
| 003 | 3 | 0 | 0 | 0 | 0 | 0 | 0 | 032 | 11 | 2 | 3 | 0 | 2 | 3 | 0 | 133 | 19 | 0 | 1 | 1 | 0 | 1 | 1 |
| 011 | 4 | 0 | 0 | 0 | 0 | 0 | 0 | 033 | 12 | 3 | 1 | 1 | 3 | 1 | 1 | 222 | 20 | 2 | 3 | 0 | 0 | 0 | 0 |
| 012 | 5 | 0 | 0 | 0 | 0 | 0 | 0 | 111 | 13 | 0 | 0 | 0 | 0 | 0 | 0 | 223 | 21 | 1 | 3 | 3 | 2 | 3 | 0 |
| 013 | 6 | 0 | 0 | 0 | 0 | 0 | 0 | 112 | 14 | 3 | 0 | 0 | 0 | 0 | 0 | 233 | 22 | 2 | 3 | 0 | 1 | 0 | 1 |
| 021 | 7 | 3 | 0 | 0 | 0 | 0 | 0 | 113 | 15 | 1 | 0 | 0 | 1 | 0 | 0 | 333 | 23 | 1 | 3 | 3 | 1 | 3 | 3 |

| rep. | k | $i_{13}$ | $i_{14}$ | $i_{15}$ | $i_{16}$ | $i_{17}$ | $i_{18}$ | rep. | k | $i_{13}$ | $i_{14}$ | $i_{15}$ | $i_{16}$ | $i_{17}$ | $i_{18}$ | rep. | k | $i_{13}$ | $i_{14}$ | $i_{15}$ | $i_{16}$ | $i_{17}$ | $i_{18}$ |
|---|---|---|---|---|---|---|---|---|---|---|---|---|---|---|---|---|---|---|---|---|---|---|---|
| 000 | 0 | 0 | 0 | 0 | 0 | 0 | 0 | 022 | 8 | 0 | 0 | 0 | 0 | 0 | 0 | 122 | 16 | 0 | 0 | 0 | 3 | 0 | 0 |
| 001 | 1 | 0 | 0 | 0 | 0 | 0 | 0 | 023 | 9 | 0 | 0 | 0 | 0 | 0 | 0 | 123 | 17 | 2 | 1 | 2 | 1 | 1 | 0 |
| 002 | 2 | 0 | 0 | 0 | 0 | 0 | 0 | 031 | 10 | 0 | 0 | 0 | 0 | 0 | 0 | 132 | 18 | 2 | 1 | 2 | 1 | 0 | 1 |
| 003 | 3 | 0 | 0 | 0 | 0 | 0 | 0 | 032 | 11 | 0 | 0 | 0 | 0 | 0 | 0 | 133 | 19 | 3 | 2 | 2 | 3 | 3 | 3 |
| 011 | 4 | 0 | 0 | 0 | 0 | 0 | 0 | 033 | 12 | 0 | 0 | 0 | 0 | 0 | 0 | 222 | 20 | 0 | 0 | 0 | 2 | 0 | 0 |
| 012 | 5 | 0 | 0 | 0 | 0 | 0 | 0 | 111 | 13 | 3 | 0 | 0 | 0 | 0 | 0 | 223 | 21 | 0 | 0 | 0 | 1 | 2 | 2 |
| 013 | 6 | 0 | 0 | 0 | 0 | 0 | 0 | 112 | 14 | 2 | 1 | 0 | 0 | 0 | 0 | 233 | 22 | 2 | 1 | 0 | 0 | 1 | 1 |
| 021 | 7 | 0 | 0 | 0 | 0 | 0 | 0 | 113 | 15 | 1 | 3 | 3 | 0 | 0 | 0 | 333 | 23 | 1 | 1 | 1 | 3 | 3 | 3 |

| rep. | k | $i_{19}$ | $i_{20}$ | $i_{21}$ | $i_{22}$ | $i_{23}$ | rep. | k | $i_{19}$ | $i_{20}$ | $i_{21}$ | $i_{22}$ | $i_{23}$ | rep. | k | $i_{19}$ | $i_{20}$ | $i_{21}$ | $i_{22}$ | $i_{23}$ |
|---|---|---|---|---|---|---|---|---|---|---|---|---|---|---|---|---|---|---|---|---|
| 000 | 0 | 0 | 0 | 0 | 0 | 0 | 022 | 8 | 0 | 0 | 0 | 0 | 0 | 122 | 16 | 0 | 0 | 0 | 0 | 0 |
| 001 | 1 | 0 | 0 | 0 | 0 | 0 | 023 | 9 | 0 | 0 | 0 | 0 | 0 | 123 | 17 | 0 | 0 | 0 | 0 | 0 |
| 002 | 2 | 0 | 0 | 0 | 0 | 0 | 031 | 10 | 0 | 0 | 0 | 0 | 0 | 132 | 18 | 0 | 0 | 0 | 0 | 0 |
| 003 | 3 | 0 | 0 | 0 | 0 | 0 | 032 | 11 | 0 | 0 | 0 | 0 | 0 | 133 | 19 | 3 | 0 | 0 | 0 | 0 |
| 011 | 4 | 0 | 0 | 0 | 0 | 0 | 033 | 12 | 0 | 0 | 0 | 0 | 0 | 222 | 20 | 0 | 1 | 0 | 0 | 0 |
| 012 | 5 | 0 | 0 | 0 | 0 | 0 | 111 | 13 | 0 | 0 | 0 | 0 | 0 | 223 | 21 | 0 | 3 | 3 | 0 | 0 |
| 013 | 6 | 0 | 0 | 0 | 0 | 0 | 112 | 14 | 0 | 0 | 0 | 0 | 0 | 233 | 22 | 2 | 1 | 2 | 1 | 0 |
| 021 | 7 | 0 | 0 | 0 | 0 | 0 | 113 | 15 | 0 | 0 | 0 | 0 | 0 | 333 | 23 | 3 | 3 | 1 | 1 | 3 |

**Example 4:**

Consider the following rotation symmetric function, where cells shaded with a same color indicate the 4 pairs of disjoint cycles where the function takes different values.

| $f$ | | | | | | | | | $x_2x_3$ | | | | | | | |
|---|---|---|---|---|---|---|---|---|---|---|---|---|---|---|---|---|
| $x_1$ | 00 | 01 | 02 | 03 | 10 | 11 | 12 | 13 | 20 | 21 | 22 | 23 | 30 | 31 | 32 | 33 |
| 0 | 0 | 1 | 2 | 3 | 1 | 0 | 1 | 2 | 2 | 2 | 0 | 1 | 3 | 1 | 2 | 0 |
| 1 | 1 | 0 | 2 | 1 | 0 | 1 | 2 | 1 | 1 | 2 | 0 | 2 | 2 | 1 | 1 | 0 |
| 2 | 2 | 1 | 0 | 2 | 2 | 2 | 0 | 1 | 0 | 0 | 2 | 1 | 1 | 2 | 1 | 0 |
| 3 | 3 | 2 | 1 | 0 | 1 | 1 | 2 | 0 | 2 | 1 | 1 | 0 | 0 | 0 | 0 | 3 |

Compact representation

| rank | 0 | 1 | 2 | 3 | 4 | 5 | 6 | 7 | 8 | 9 | 10 | 11 | 12 | 13 | 14 | 15 | 16 | 17 | 18 | 19 | 20 | 21 | 22 | 23 |
|---|---|---|---|---|---|---|---|---|---|---|---|---|---|---|---|---|---|---|---|---|---|---|---|---|
| $\langle \mathbf{F} \rangle$ | 0 | 1 | 2 | 3 | 0 | 1 | 2 | 2 | 0 | 1 | 1 | 2 | 0 | 1 | 2 | 1 | 0 | 2 | 1 | 0 | 2 | 1 | 0 | 3 |

RMF spectrum of $f$

| $S_f$ | | | | | | | | $x_2x_3$ | | | | | | | | |
|---|---|---|---|---|---|---|---|---|---|---|---|---|---|---|---|---|
| $x_1$ | 00 | 01 | 02 | 03 | 10 | 11 | 12 | 13 | 20 | 21 | 22 | 23 | 30 | 31 | 32 | 33 |
| 0 | 0 | 3 | 0 | 0 | 3 | 2 | 2 | 2 | 0 | 1 | 2 | 3 | 0 | 2 | 1 | 1 |
| 1 | 3 | 2 | 1 | 2 | 2 | 0 | 3 | 2 | 2 | 3 | 2 | 0 | 2 | 2 | 0 | 0 |
| 2 | 0 | 2 | 2 | 1 | 1 | 3 | 2 | 0 | 2 | 2 | 0 | 1 | 3 | 0 | 1 | 0 |
| 3 | 0 | 2 | 3 | 1 | 2 | 2 | 0 | 0 | 1 | 0 | 1 | 0 | 1 | 0 | 0 | 0 |

(Notice that $S_f$ is "less rotation symmetric" than $f$ in the sense that in the blocks with rank 6 and 10 and with rank 17 and 18, the spectrum does not have respectively different values)

Compact representation

| rank | 0 | 1 | 2 | 3 | 4 | 5 | 6 | 7 | 8 | 9 | 10 | 11 | 12 | 13 | 14 | 15 | 16 | 17 | 18 | 19 | 20 | 21 | 22 | 23 |
|---|---|---|---|---|---|---|---|---|---|---|---|---|---|---|---|---|---|---|---|---|---|---|---|---|
| $\langle S_f \rangle$ | 0 | 3 | 0 | 0 | 2 | 2 | 2 | 1 | 2 | 3 | 2 | 1 | 1 | 0 | 3 | 2 | 2 | 0 | 0 | 0 | 0 | 1 | 0 | 0 |

From Theorem 6:

$\langle S_f \rangle = \langle SB(4,1) \rangle \oplus 2\langle SB(4,2) \rangle \oplus 3\langle SB(4,3) \rangle \oplus \langle SB(4,5) \rangle \oplus 2\langle SB(4,6) \rangle \oplus 2\langle SB(4,7) \rangle \oplus \langle SB(4,9) \rangle \oplus \langle SB(4,10) \rangle \oplus$
$\oplus 2\langle SB(4,11) \rangle \oplus \langle SB(4,13) \rangle \oplus 2\langle SB(4,14) \rangle \oplus \langle SB(4,15) \rangle \oplus 2\langle SB(4,17) \rangle \oplus \langle SB(4,18) \rangle \oplus 2\langle SB(4,20) \rangle \oplus$
$\oplus \langle SB(4,21) \rangle \oplus 3\langle SB(4,23) \rangle \mod 4$

From Corollary 6.1:

$\langle \mathbf{F} \rangle = 3\langle SB(4,1) \rangle \oplus 2\langle SB(4,4) \rangle \oplus 2\langle SB(4,5) \rangle \oplus 2\langle SB(4,6) \rangle \oplus \langle SB(4,7) \rangle \oplus 2\langle SB(4,8) \rangle \oplus 3\langle SB(4,9) \rangle \oplus 2\langle SB(4,10) \rangle \oplus$
$\oplus \langle SB(4,11) \rangle \oplus \langle SB(4,12) \rangle \oplus 3\langle SB(4,14) \rangle \oplus 2\langle SB(4,15) \rangle \oplus 2\langle SB(4,16) \rangle \oplus \langle SB(4,21) \rangle \mod 4$

Both expressions give the correct result.

**Example 5:**

A function was built based on the function of Example 4, but in each pair of disjoint cycles, in one cycle the value was changed. In the "only symmetric" part of the function at a few random places the value was changed, but preserving the symmetry.

| $f$ | | | | | | | | $x_2x_3$ | | | | | | | | |
|---|---|---|---|---|---|---|---|---|---|---|---|---|---|---|---|---|
| $x_1$ | 00 | 01 | 02 | 03 | 10 | 11 | 12 | 13 | 20 | 21 | 22 | 23 | 30 | 31 | 32 | 33 |
| 0 | 0 | 3 | 2 | 3 | 3 | 1 | 3 | 0 | 2 | 2 | 1 | 0 | 3 | 1 | 2 | 1 |
| 1 | 3 | 1 | 2 | 1 | 1 | 1 | 2 | 1 | 3 | 2 | 0 | 3 | 0 | 1 | 1 | 0 |
| 2 | 2 | 3 | 1 | 2 | 2 | 2 | 0 | 1 | 1 | 0 | 2 | 1 | 0 | 3 | 1 | 0 |
| 3 | 3 | 0 | 0 | 1 | 1 | 1 | 3 | 0 | 2 | 1 | 1 | 0 | 1 | 0 | 0 | 3 |

Compact representation:

| rank | 0 | 1 | 2 | 3 | 4 | 5 | 6 | 7 | 8 | 9 | 10 | 11 | 12 | 13 | 14 | 15 | 16 | 17 | 18 | 19 | 20 | 21 | 22 | 23 |
|---|---|---|---|---|---|---|---|---|---|---|---|---|---|---|---|---|---|---|---|---|---|---|---|---|
| $\langle \mathbf{F} \rangle$ | 0 | 3 | 2 | 3 | 1 | 3 | 0 | 2 | 1 | 0 | 1 | 2 | 1 | 1 | 2 | 1 | 0 | 3 | 1 | 0 | 2 | 1 | 0 | 3 |

RMF spectrum

| $S_f$ | | | | | | | | | $x_2x_3$ | | | | | | | |
|---|---|---|---|---|---|---|---|---|---|---|---|---|---|---|---|---|
| $x_1$ | 00 | 01 | 02 | 03 | 10 | 11 | 12 | 13 | 20 | 21 | 22 | 23 | 30 | 31 | 32 | 33 |
| 0 | 0 | 1 | 0 | 2 | 1 | 3 | 0 | 1 | 0 | 1 | 3 | 3 | 2 | 1 | 2 | 3 |
| 1 | 1 | 3 | 1 | 1 | 3 | 1 | 2 | 3 | 0 | 2 | 3 | 2 | 1 | 3 | 0 | 1 |
| 2 | 0 | 0 | 3 | 2 | 1 | 2 | 3 | 0 | 3 | 3 | 3 | 3 | 3 | 2 | 3 | 0 |
| 3 | 2 | 1 | 3 | 3 | 1 | 3 | 2 | 1 | 2 | 0 | 3 | 0 | 3 | 1 | 0 | 3 |

Compact representation

| rank | 0 | 1 | 2 | 3 | 4 | 5 | 6 | 7 | 8 | 9 | 10 | 11 | 12 | 13 | 14 | 15 | 16 | 17 | 18 | 19 | 20 | 21 | 22 | 23 |
|---|---|---|---|---|---|---|---|---|---|---|---|---|---|---|---|---|---|---|---|---|---|---|---|---|
| $\langle S_f \rangle$ | 0 | 1 | 0 | 2 | 3 | 0 | 1 | 1 | 3 | 3 | 1 | 2 | 3 | 1 | 2 | 3 | 3 | 2 | 0 | 1 | 3 | 3 | 0 | 3 |

The function and the RMF spectrum satisfy Theorem 6 and Corollary 6.1

$$\langle \mathbf{F} \rangle = \bigoplus_{k=0}^{k=23} \langle S_f \rangle_k \langle SB(4,k) \rangle \mod 4 \qquad \langle S_f \rangle = \bigoplus_{k=0}^{k=23} \langle \mathbf{F} \rangle_k \langle SB(4,k) \rangle \mod 4$$

**Example 6:**

| $f$ | | | | | | | | | $x_2x_3$ | | | | | | | |
|---|---|---|---|---|---|---|---|---|---|---|---|---|---|---|---|---|
| $x_1$ | 00 | 01 | 02 | 03 | 10 | 11 | 12 | 13 | 20 | 21 | 22 | 23 | 30 | 31 | 32 | 33 |
| 0 | 0 | 1 | 2 | 3 | 1 | 0 | 1 | 2 | 2 | 2 | 0 | 1 | 3 | 1 | 2 | 0 |
| 1 | 1 | 0 | 2 | 1 | 0 | 1 | 2 | 1 | 1 | 2 | 1 | 3 | 2 | 1 | 1 | 0 |
| 2 | 2 | 1 | 0 | 2 | 2 | 2 | 1 | 1 | 0 | 1 | 2 | 1 | 1 | 3 | 1 | 0 |
| 3 | 3 | 2 | 1 | 0 | 1 | 1 | 3 | 0 | 2 | 1 | 1 | 0 | 0 | 0 | 0 | 3 |

Compact representation:

| rank | 0 | 1 | 2 | 3 | 4 | 5 | 6 | 7 | 8 | 9 | 10 | 11 | 12 | 13 | 14 | 15 | 16 | 17 | 18 | 19 | 20 | 21 | 22 | 23 |
|---|---|---|---|---|---|---|---|---|---|---|---|---|---|---|---|---|---|---|---|---|---|---|---|---|
| $\langle \mathbf{F} \rangle$ | 0 | 1 | 2 | 3 | 0 | 1 | 2 | 2 | 0 | 1 | 1 | 2 | 0 | 1 | 2 | 1 | 1 | 3 | 1 | 0 | 2 | 1 | 0 | 3 |

RMF spectrum;

| $S_f$ | | | | | | | | | $x_2x_3$ | | | | | | | |
|---|---|---|---|---|---|---|---|---|---|---|---|---|---|---|---|---|
| $x_1$ | 00 | 01 | 02 | 03 | 10 | 11 | 12 | 13 | 20 | 21 | 22 | 23 | 30 | 31 | 32 | 33 |
| 0 | 0 | 3 | 0 | 0 | 3 | 2 | 2 | 2 | 0 | 1 | 2 | 3 | 0 | 2 | 1 | 1 |
| 1 | 3 | 2 | 1 | 2 | 2 | 0 | 3 | 2 | 2 | 3 | 1 | 2 | 2 | 2 | 1 | 2 |
| 2 | 0 | 2 | 2 | 1 | 1 | 3 | 1 | 1 | 2 | 1 | 2 | 0 | 3 | 2 | 0 | 1 |
| 3 | 0 | 2 | 3 | 1 | 2 | 2 | 2 | 2 | 1 | 1 | 0 | 1 | 1 | 2 | 1 | 2 |

Compact representation:

| rank | 0 | 1 | 2 | 3 | 4 | 5 | 6 | 7 | 8 | 9 | 10 | 11 | 12 | 13 | 14 | 15 | 16 | 17 | 18 | 19 | 20 | 21 | 22 | 23 |
|---|---|---|---|---|---|---|---|---|---|---|---|---|---|---|---|---|---|---|---|---|---|---|---|---|

| ⟨$S_f$⟩ | 0 | 3 | 0 | 0 | 2 | 2 | 2 | 1 | 2 | 3 | 2 | 1 | 1 | 0 | 3 | 2 | 1 | 2 | 1 | 2 | 2 | 0 | 1 | 2 |
|---|---|---|---|---|---|---|---|---|---|---|---|---|---|---|---|---|---|---|---|---|---|---|---|---|

Notice that the spectrum has the value "2" in roughly 40% of the places, meaning that when evaluating the equation of Corollary 6.1, a considerable amount of $SB(4,k)$ coefficients will be scaled by 2; and 2 is a zero divider in $(\mathbb{Z}_4, \oplus, \cdot)$. The corollary was satisfied. It may be concluded that the zero divider does not affect the Theorem and its Corollary.

**Example 7:** The following function is rotation symmetric and except for the necessary disjoint cycles with different values, and the places where the cycle has size 1, in all other places the function has the value 2.

| $f$ | | | | | | | | $x_2x_3$ | | | | | | | | |
|---|---|---|---|---|---|---|---|---|---|---|---|---|---|---|---|---|
| $x_1$ | 00 | 01 | 02 | 03 | 10 | 11 | 12 | 13 | 20 | 21 | 22 | 23 | 30 | 31 | 32 | 33 |
| 0 | 0 | 2 | 2 | 2 | 2 | 2 | 1 | 2 | 2 | 2 | 2 | 1 | 2 | 1 | 2 | 2 |
| 1 | 2 | 2 | 2 | 1 | 2 | 1 | 2 | 2 | 1 | 2 | 2 | 2 | 2 | 2 | 1 | 2 |
| 2 | 2 | 1 | 2 | 2 | 2 | 2 | 2 | 1 | 2 | 2 | 2 | 2 | 1 | 2 | 2 | 2 |
| 3 | 2 | 2 | 1 | 2 | 1 | 2 | 2 | 2 | 2 | 1 | 2 | 2 | 2 | 2 | 2 | 3 |

Compact representation

| rank | 0 | 1 | 2 | 3 | 4 | 5 | 6 | 7 | 8 | 9 | 10 | 11 | 12 | 13 | 14 | 15 | 16 | 17 | 18 | 19 | 20 | 21 | 22 | 23 |
|---|---|---|---|---|---|---|---|---|---|---|---|---|---|---|---|---|---|---|---|---|---|---|---|---|
| ⟨**F**⟩ | 0 | 2 | 2 | 2 | 2 | 1 | 2 | 2 | 2 | 1 | 1 | 2 | 2 | 1 | 2 | 2 | 2 | 2 | 1 | 2 | 2 | 2 | 2 | 3 |

RMF Spectrum

| $S_f$ | | | | | | | | $x_2x_3$ | | | | | | | | |
|---|---|---|---|---|---|---|---|---|---|---|---|---|---|---|---|---|
| $x_1$ | 00 | 01 | 02 | 03 | 10 | 11 | 12 | 13 | 20 | 21 | 22 | 23 | 30 | 31 | 32 | 33 |
| 0 | 0 | 2 | 2 | 2 | 2 | 2 | 3 | 1 | 2 | 2 | 0 | 1 | 2 | 1 | 3 | 3 |
| 1 | 2 | 2 | 2 | 1 | 2 | 3 | 1 | 3 | 3 | 1 | 1 | 3 | 1 | 3 | 3 | 3 |
| 2 | 2 | 3 | 0 | 3 | 2 | 1 | 1 | 3 | 0 | 1 | 0 | 2 | 1 | 3 | 2 | 0 |
| 3 | 2 | 1 | 1 | 3 | 1 | 3 | 3 | 3 | 3 | 3 | 2 | 0 | 3 | 3 | 0 | 0 |

The RMF spectrum is rotation symmetric, although only in two pairs of disjoint cycles it takes different values.

Compact representation

| rank | 0 | 1 | 2 | 3 | 4 | 5 | 6 | 7 | 8 | 9 | 10 | 11 | 12 | 13 | 14 | 15 | 16 | 17 | 18 | 19 | 20 | 21 | 22 | 23 |
|---|---|---|---|---|---|---|---|---|---|---|---|---|---|---|---|---|---|---|---|---|---|---|---|---|
| ⟨$S_f$⟩ | 0 | 2 | 2 | 2 | 2 | 3 | 1 | 2 | 0 | 1 | 1 | 3 | 3 | 3 | 1 | 3 | 1 | 3 | 3 | 3 | 0 | 2 | 0 | 0 |

Both weighted sums required by Theorem 6 and Corollary 6.1 give the correct results. The zero divider "2" does not affect Theorem 6.